\documentclass[12pt]{article}
\RequirePackage{amsthm,amsmath}
\usepackage{amsfonts}
\usepackage{amssymb}
\usepackage{fancyvrb}
\usepackage{graphicx}
\usepackage{url}
\DeclareSymbolFont{AMSb}{U}{msb}{m}{n}
\DeclareMathSymbol{\N}{\mathbin}{AMSb}{"4E}
\DeclareMathSymbol{\Z}{\mathbin}{AMSb}{"5A}
\DeclareMathSymbol{\R}{\mathbin}{AMSb}{"52}
\DeclareMathSymbol{\Q}{\mathbin}{AMSb}{"51}
\DeclareMathSymbol{\I}{\mathbin}{AMSb}{"49}
\DeclareMathSymbol{\C}{\mathbin}{AMSb}{"43}

\title{A Distributed Procedure for Computing Stochastic Expansions with \texttt{Mathematica}}

\author{ Christophe Ladroue\footnote{Department of Statistics, University of Warwick, Coventry, CV4 7AL, UK.}\\ and \\ Anastasia Papavasiliou\footnote{Department of Statistics, University of Warwick, Coventry, CV4 7AL, UK.}}

\begin{document}
\maketitle
\abstract{
The solution of a (stochastic) differential equation can be locally approximated by a (stochastic) expansion. 
If the vector field of the differential equation is a polynomial, the corresponding expansion is a linear combination of iterated integrals of the drivers and can be calculated using Picard Iterations.
However, such expansions grow exponentially fast in their number of terms, due to their specific algebra, rendering their practical use limited. 

We present a \texttt{Mathematica} procedure that addresses this issue by re-parametrising the polynomials and distributing the load in as small as possible parts that can be processed and manipulated independently, thus alleviating large memory requirements and being perfectly suited for parallelized computation. We also present an iterative implementation of the shuffle product (as opposed to a recursive one, more usually implemented) as well as a fast way for calculating the expectation of iterated Stratonovich integrals for Brownian Motion.
}



					\section{Motivation and mathematical background}
In this section, we introduce the mathematical background and motivation for manipulating expansions and iterated integrals. 
The next subsection introduces the Picard procedure, a simple iterative way to derive local approximation of the solution of a differential equation. Iterated integrals are then introduced and the two are finally combined to define the expansions.

					\subsection{Motivation and notation}
Consider the following differential equation:
\begin{equation}\label{eq:mainSDE}
dY_t=f(Y_t,\theta)dX_t
\end{equation}

where $Y_t\in\R^m$, $X_t\in\R^n$ and $f(.,\theta)\in\R^{m\times n}$. The parameters of the function $f$ are collected in the variable $\theta$. We call each $X^{(i)}$ a driver of the differential equation. We assume the functions $f_{j,i},\ i\in \{1\hdots n\},\ j\in\{1\hdots m\}$ to be polynomials. The initial value of $Y_t$ is set to $Y_0$.

The objective is to derive a local approximation of the solution in terms of $f$, $Y_0$ and the iterated integrals of $Y_t$, \emph{i.e.}that is integrals of the form 
\[\int\hdots \int_{0< u_1< \cdots < u_k \leq T} dX^{\tau_1}_{u_1}\hdots dX^{\tau_k}_{u_k}\]
 for any word $\tau=(\tau_1,\dots,\tau_k)$ constituted of the letter $\tau_i\in\{1,\hdots,n\},\ i = 1,\dots,k.$. Such expansions play an important role in the theory of rough paths, allowing one to define such a differential equation for a large class of drivers $X$ (\cite{Lyons2003}). They have also been used recently for parameter estimation of SDE (\cite{Papavasiliou2010}).

					\subsection{Picard iterations}
Picard iterations provide a way for deriving local approximations of solutions of differential equation. They are defined as:
\begin{eqnarray}\label{eq:picard}
Y_{0,T}(0)&=&Y_0-Y_0=0\\
Y_{0,T}(r+1,j)&=&\sum_{i=1}^{n}\int f_{j,i}(Y_s(r,j))dX_s^{(i)}\\ 
\end{eqnarray}
where $Y_{0,T}(r,j)$ is the $j^{th}$ component of the approximation of $Y_{0,T}=\int_0^TdY_s=Y_T-Y_0$ after $r$ iterations. Thus, the first iteration gives:
\begin{equation}\label{eq:picardFirstIteration}
Y_{0,T}(1,j)=\int_0^t f_{j,1}(Y_0)dX_s^{(1)}+\hdots+\int_0^t f_{j,n}(Y_0)dX_s^{(n)}
\end{equation}
Note that if we are interested in the actual value of the solution at a time $t\in [0,T]$, its Picard approximation is given by
\[ Y_t(r,j) = Y_0 + Y_{0,t}(r,j).\]
One of the main successes in the theory of rough paths was to give precise conditions on $X$ and $f$ for Picard iterations to converge (\cite{Lyons2003}).

					\subsection{Iterated integrals}
Iterated integrals are integrals of the form:
$$X_{s,t}^{(\tau)}=\int \hdots \int_{s< u_1< \hdots < u_k<t}dX_{u_1}^{(\tau_1)}\hdots dX_{u_k}^{(\tau_k)}$$
where $\tau=(\tau_1,\dots,\tau_k)$ is called a word, with letters $\tau_i\in\{1\hdots n\}, i = 1,\dots,k$. By definition, integrating an iterated integral produces an iterated integral:
\begin{eqnarray}
\int_0^T X_{0,s}^{(\tau)}dX_s^{(j)}&=&\int_0^T \int \hdots \int_{0< u_1< \hdots < u_k<s}dX_{u_1}^{(\tau_1)}\hdots dX_{u_k}^{(\tau_k)} dX_s^{(j)}\\
&=&X_{0,T}^{(\tau_1,\hdots,\tau_k, j)}
\end{eqnarray}
If $X$ is a geometric $p$-rough path, \emph{i.e.} it can be approximated by paths of bounded variation (see \cite{Lyons2003} for a precise definition), then the integrals obey the usual integration-by-parts rule, which can be generalized as follows:
\begin{equation}\label{eq:shuffle}
X_{s,t}^{(\tau)} X_{s,t}^{(\rho)}=\sum_{\alpha\in\tau\sqcup\rho} X_{s,t}^{(\alpha)}
\end{equation}
The shuffle product $\sqcup$ of two words $\tau$ and $\rho$ is the set of all words using the letters in $\tau$ and $\rho$ such that they are in their original
order. For example $13\sqcup 42 = \{1342, 1432, 4132, 1423, 4213,4132\}$ but 1324, for example, does not belong to the set. 

This result holds only in the case of deterministic drivers X, or Stratonovitch integrals. A similar relation exists for It\^o integrals but requires a small correcting term (\cite{Tocino2009}) making them slightly less practical in this context. A simple transformation allows diffusion processes to be defined with respect to It\^o or Stratonovitch integrals (\cite{Kloeden91}):
$$
dy_t=\mu(y_t)dt+\sigma(y_t) dB_t \Leftrightarrow dy_t=\left(\mu(y_t)-\frac{1}{2}\sigma'(y_t)\sigma(y_t)\right)dt+\sigma(y_t)\circ dB_t 
$$

It immediately follows from equation (\ref{eq:shuffle}) that any power of an iterated integral is a sum of iterated integrals and that a polynomial of iterated integrals is also a linear combination of single iterated integrals. For example:
\begin{eqnarray*}
X_{0,T}^{(1)}X_{0,T}^{(2,3)}&=&X_{0,T}^{(1,2,3)}+X_{0,T}^{(2,1,3)}+X_{0,T}^{(2,3,1)}\\
\left(X_{0,T}^{(1)}\right)^2&=&X_{0,T}^{(1)} X_{0,T}^{(1)}\\
&=&X_{0,T}^{(1,1)}+X_{0,T}^{(1,1)}\\
&=&2X_{0,T}^{(1,1)}\\
\left(X_{0,T}^{(1)}\right)^\ell&=&\ell!\ X_{0,T}^{(1,\hdots,1)}\\
X_{0,T}^{(0,1,0)}X_{0,T}^{(1,1)}&=&X_{0,T}^{(0, 1, 0, 1, 1)} + 2 X_{0,T}^{(0, 1, 1, 0, 1)} + 3 X_{0,T}^{(0, 1, 1, 1, 0)} \\
&& + X_{0,T}^{(1, 0, 1, 0, 1)} + 2 X_{0,T}^{(1, 0, 1, 1, 0)} + X_{0,T}^{(1, 1, 0, 1, 0)}\\
\end{eqnarray*}
The number of terms from the product of iterated integrals grows fast, exponentially if all letters are different.
					\subsection{Expansions}
We are now in position to derive expansions for the solution of a differential equation. Picard iterations yield:
\begin{eqnarray*}
Y_{0,T}(0,j)&=&0\\
Y_{0,T}(1,j)&=&\int_{0,T}f_{j1}(Y_s(0))dX_s^{(1)}+\hdots+\int_{0,T}f_{jn}(Y_s(0))dX_s^{(n)}\\
&=&f_{j1}(Y_0)X_{0,T}^{(1)}+\hdots+f_{jn}(Y_0)X_{0,T}^{(n)}\\
Y_{0,T}(2,j)&=&\int_{0,T}f_{j1}(Y_s(1))dX_s^{(1)}+\hdots+\int_{0,T}f_{jn}(Y_s(1))dX_s^{(n)}\\
&=&\int_{0,T}f_{j1}(Y_{0,s}(1)+Y_0)dX_s^{(1)}+\hdots+\int_{0,T}f_{jn}(Y_s(1)+Y_0)dX_s^{(n)}
\end{eqnarray*}
Since $Y_{0,s}(1,.)$ is a sum of iterated integrals, the polynomials $f_{ji}(Y_s(1)+Y_0)$ are also sums of iterated integrals and so is their integration with respect to $X^{(i)}$. $Y_{0,T}(2)$ is thus a sum of iterated integrals and by recursion all $Y_{0,T}(r,.)$ are.  A formal proof in the context of differential equations driven by rough paths can be found in \cite{Papavasiliou2010}.

					\subsubsection*{Example}
Consider the Ornstein-Uhlenbeck process: $dy_t=a(1-y_t)dt+bdW_t$ and $y_0=0$. In this case, $X_t^{(1)}=t$, $X_t^{(2)}=W_t$ and $Y_0^{(i)}=0$ for $i\in\{1,2\}$. Applying Picard iterations, we obtain:
\begin{eqnarray*}
Y_{0,T}(0)&=&0\\
Y_{0,T}(1)&=&\int_0^Ta(1-0)dX_s^{(1)}+\int_0^TbdX_s^{(2)}\\
&=& aX_{0,T}^{(1)}+bX_{0,T}^{(2)}\\
Y_{0,T}(2)&=&\int_0^Ta(1-(aX_s^{(1)}+bX_s^{(2)}))dX_s^{(1)}+\int_0^TbdX_s^{(2)}\\
&=&aX_{0,T}^{(1)}-a^2X_{0,T}^{(1,1)}-abX_{0,T}^{(2,1)}+bX_{0,T}^{(2)}\\
Y_{0,T}(3)&=&aX_{0,T}^{(1)}-a^2X_{0,T}^{(1,1)}+a^3X_{0,T}^{(1,1,1)}+a^bX_{0,T}^{(2,1,1)}+abX_{0,T}^{(2,1)}+bX_{0,T}^{(2)}
\end{eqnarray*}
The solution of the stochastic differential equation can thus be approximated by a series of iterated integral of the drivers, whose coefficients are a function of the parameters. The iterated integrals capture the statistics of the drivers and are separated from the parameters. 

This derivation can be readily implemented in \texttt{Mathematica} (\cite{Wolfram2003}\footnote{The code will work with \texttt{Mathematica} version 7, when parallel computing was introduced.}) (see \cite{Tocino2009} for an implementation of the shuffle product) but suffers a major drawback: each product of iterated integrals being a shuffle product, the number of terms produced grows extremely fast (exponentially in the worst cases) and rapidly becomes unmanageable. In the next section, we introduce a re-parametrisation of the problem that circumvents this problem by providing an alternative representation of the expansion which can be processed in a distributed manner, alleviating large memory requirements.

\section{Re-parametrisation of the polynomials}
One thing to note in the derivation of expansions is that each successive iterations requires the explicit linear combination of iterated integrals for the previous iteration; evaluating $Y_{0,T}(r)$ requires the complete expansion for $Y_{0,T}(r-1)$. As the number of terms grows, manipulating this object rapidly becomes unwieldy.

In this section, we describe a re-parametrisation of the polynomials which bypasses the need for an expansion in terms of iterated integrals. It provides a more compact representation of the approximate solution and is naturally amenable to parallel processing. For clarity of exposition, we first introduce the approach in the case of a one-dimensional ($m=1$) differential equation with $n$ drivers. We assume the polynomials $f_{1i}$ to be of degree less or equal to $q$. In the last subsection, the procedure is generalized to $m$-dimensional differential equations.

					\subsection{One-dimensional case}
We first remark that a polynomial $P(y)$ can be written in terms of $y-y_0$ by writing its Taylor expansion around $y_0$:
\begin{equation}\label{eq:polynomialTaylor}
P(y)=\sum_{k=0}^q\frac{1}{k!}\partial^kP(y_0)(y-y_0)^k
\end{equation}
Next we introduce the following new operation for iterated integrals:
\begin{equation}\label{eq:productDefinition}
X_{s,t}^\alpha\vartriangleright X_{s,t}^\beta=\int_s^tX_{s,u}^\alpha dX_{s,u}^{\beta}=\int_s^tX_{s,u}^\alpha X_{s,u}^{\beta_-}dX_u^{\beta_{\textnormal{end}}}
\end{equation}
where $\beta_-$ is the word $\beta$ with the last letter removed and $\beta_\textnormal{end}$ the last letter of $\beta$. This is a non-associative, non-commutative operation -- in fact, it can be viewed as a non-commutative dendriform. We can rewrite the Picard iteration using this operation. Note that from now on, the interval $[s,t]$ will be fixed to $[0,T]$ and will be omitted.
\begin{eqnarray}
Y(0)&=&y_0-y_0=0\\
Y(1)&=&\sum_{i=1}^n\int_0^Tf_{j1}(Y_s(0)+y_0)dX_s^{(i)}\\
&=&\sum_{i=1}^nf_{j1}(y_0)X^{(i)}\\
Y(r+1)&=&\sum_{i=1}^nf_{1i}(Y_s(r)+y_0)\vartriangleright X^{(i)}\\
&=&\sum_{i=1}^n \sum_{k=0}^q\frac{1}{k!}\partial^kf_{1i}(y_0)Y_s(r)^k\vartriangleright X^{(i)}\\
&=&\sum_{k=0}^q \left( ( (Y_s(r))^k\vartriangleright \sum_{i=1}^n (\partial^kf_{1i}(y_0)X^{(i)})\right) 
\end{eqnarray}

Therefore, if we define the objects $Q$ as:
\begin{equation}\label{eq:QDefinition}
Q^k=\sum_{i=1}^n\frac{1}{k!}\partial^kf_{1i}(y_0)X^{(i)},
\end{equation}
the Picard iteration takes the following form:
\begin{eqnarray*}
Y(1)&=&Q^0\\
Y(r+1)&=&\sum_{k=0}^qY(r)^k\vartriangleright Q^k\\
\end{eqnarray*}
where $Y(r)^k$ is the usual product.
					\subsection{Description of the approach}\label{subsection:approach}
Consider a one-dimensional differential equation with quadratic functions $f_{1i}$, \emph{i.e.} $q=2$. The new representation yields:
\begin{eqnarray}
Y(1)&=&Q^0\\
Y(2)&=&1\vartriangleright Q^0+Q^0\vartriangleright Q^1+(Q^0)^2\vartriangleright Q^2\\
&=&Q^0+Q^0\vartriangleright Q^1+(Q^0)^2\vartriangleright Q^2\\
Y(3)&=&Q^0+(Q^0+Q^0\vartriangleright Q^1+(Q^0)^2\vartriangleright Q^2)\vartriangleright Q^1\\
&&+(Q^0+Q^0\vartriangleright Q^1+(Q^0)^2\vartriangleright Q^2)^2\vartriangleright Q^2
\end{eqnarray}
where $(Q^k)^\ell$ is the $k^{th}$ object $Q$ to the power $\ell$. Crucially, this new representation does not require the explicit computation of the shuffle products, keeping the number of terms under control. Moreover, the expression can be expanded into its summands and each summand be processed independently. Thus, we avoid the handling of a large expression and are able to parallelize the computation. 

Note that the representation only depends on the maximum degree of the polynomials but not on the number of drivers. The Picard iteration can therefore be done only once, stored in a file and used at a later date for any system that uses the same maximum degree $q$. 

Using this representation, it is possible to derive the stochastic expansions through a few stages:
\begin{enumerate}
\item Expand the expression $Y(r)$ into its mononials $u$. Each monomial $u$ is a function of $Q$'s that uses the non-commutative products. Importantly, the objects $Q$ and the product $\vartriangleright$ are only used as place-holders at that stage. For example, the first three monomials for $Y(3)$ are $Q^0$, $Q^0\vartriangleright Q^1$ and $(Q^0\vartriangleright Q^1)\vartriangleright Q^1$.
\item For each monomial $u$, the objects $Q$ are replaced by their values from the model at hand. Each $Q$ is a weighted sum of the drivers $X^{(i)}$ (eq.\ref{eq:QDefinition}), so each $u$ becomes a polynomial $V$ of $X^{(i)}$ in terms of the non-commutative product. As in the previous step, the product is still only employed as a place-holder and not instantiated.
\item Each polynomial $V$ is expanded into its monomials $v$. Each $v$ is a function of $X^{(i)}$ and $\vartriangleright$.
\item For each monomial $v$, the product is instantiated; its actual definition in terms of shuffle product is only used at this later stage.
\end{enumerate}

Each process only requires a fraction of the memory that a direct approach (replacing $Q$'s by the $X^{(i)}$ and using the product's definition) would. Moreover, at all stages, each monomial can be processed independently from the rest, leading to natural parallelization. It is also important to note that the whole expression is actually stored in a file and thus not in memory at any point. The exact details of this procedure are given in section \ref{sec:implementation}.

					\subsection{Generalization}\label{section:generalization}
In the multidimensional case, the Taylor expansion of a polynomial requires a larger number of terms that involve cross-products between the components of the vector $Y$. The objects $Q$ are not indexed by $k\in\{0,\hdots,q\}$ but by the set $\textnormal{OW}_m(0,q)$ of the ordered words of length up to $q$ written with letters $\{1,\hdots,m\}$ and are now defined as:
\begin{equation}
Q^\tau_j=\sum_{i=1}^n\frac{|\tau|!}{c(\tau)}\partial_\tau f_{ji}(Y_0)X^{(i)}
\end{equation} 
for $j\in\{1,\hdots,m\}$. The constant $c(\tau)$ is the number of different words we can construct using the letters in $\tau$. The Picard iteration becomes:
$$
Y(1)^{(j)}=Q^\emptyset_j\ \textnormal{and}\ Y(r+1)^{(j)}=\sum_{\tau\in\textnormal{OW}_m(0,q)}Y(r)^\tau\vartriangleright Q^\tau_j
$$

\section{Implementation}\label{sec:implementation}
This section describes how this new approach was implemented in \texttt{Mathematica}. 
In this part, iterated integrals of the drivers ($X^{(i_1,\hdots,i_n)}$) are denoted $j^{(i_1,\hdots,i_n)}$ to follow convention and drivers are numbered from $0$ with the first driver representing time.  
					\subsection{Shuffle product}
Since each product of two iterated integrals is a shuffle product, special care must be taken of its implementation. \cite{Tocino2009} has shown a way of writing the product in \texttt{Mathematica}:
\begin{Verbatim}
Tocino[j[{x_}], j[a_List]] := 
  Sum[j[Insert[a, x, k]], {k, 1, Length@a + 1}];
Tocino[j[a_List], j[{x_}]] := Tocino[j[{x}], j[a]];
Tocino[j[a_List], j[b_List]] :=
    Ap[Tocino[j[a], j[Drop[b, -1]]], Last@b] + 
	    Ap[Tocino[j[Drop[a, -1]], j[b]], Last@a] /; (Length@a > 1 && Length@b > 1);
\end{Verbatim}
This is a direct and natural translation of the following result:
$$
J^\alpha J^\beta=\int J^{\alpha-}J^\beta dJ^{\alpha_\textnormal{end}}+\int J^\alpha J^{\beta-} dJ^{\beta_\textnormal{end}}
$$

We present a new implementation of the shuffle product. The product is done iteratively instead of recursively and is based on string transformation. To calculate the shuffle product of two words $\alpha$ and $\beta$, we first set the string `aa...aabb..bb', with as many $a$'s and $b$'s as there are letters in $\alpha$ and $\beta$ respectively (line 2). We then replace all occurences of `ab' by `ba' (line 6) and keep on iterating the replacement until none are possible (\emph{i.e} when the word is `bb..bbaa..aa') while keeping track of the new words generated (lines 4-6). The set of all words thus created is the shuffle product of two arbitrary words of the same lengths of $\alpha$ and $\beta$. Finally, the letters $a$ and $b$ in each word are replaced by their actual values from $\alpha$ and $\beta$ (lines 7-9).

\begin{Verbatim}[numbers=left]
Shuffle[a_List,b_List]:=Module[{list,u,v},
	u={StringJoin@Join[Table["a",{Length[a]}],Table["b",{Length[b]}]]};
	v=Table[0,{StringLength[u[[1]]]}];
	list=Flatten@NestWhileList[
			DeleteDuplicates@Flatten@
				StringReplaceList[#,"ab"->"ba"]&,u,#!={}&];
	(v[[#[[1]]&/@StringPosition[#,"a"]]]=a;
	v[[#[[1]]&/@StringPosition[#,"b"]]]=b;
	v)&/@list];
\end{Verbatim}

\begin{figure}
\label{fig:shuffles}
\begin{center}
\includegraphics[width=0.75\textwidth]{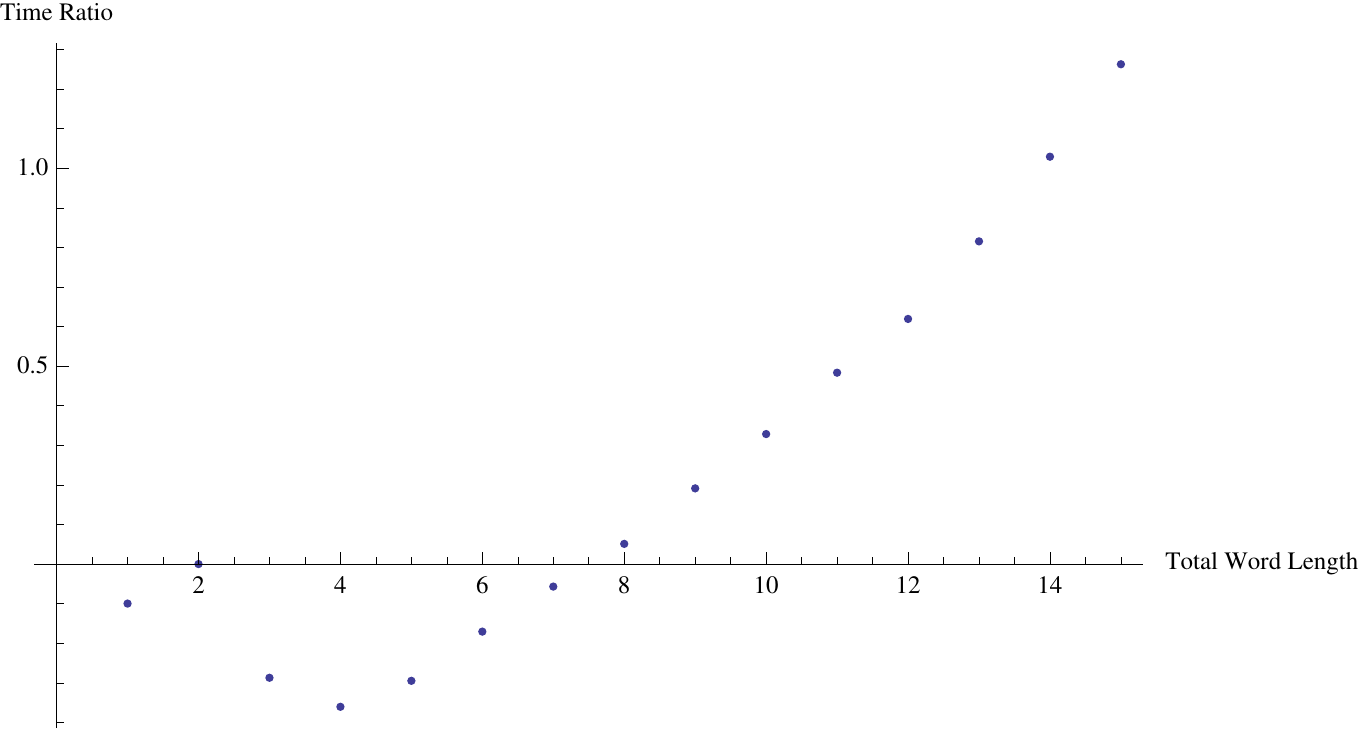}
\end{center}
\caption{Log2 of the average of the relative speed of the iterative over the recursive method.}
\end{figure}

Figure \ref{fig:shuffles} shows the relative performances of the implementation. The product of two iterated integrals with random words is calculated with the recursive definition and the iterative definition. The ratio of the time spent by the former over the latter is recorded and the process is re-iterated 1,000 times for each final-word length. On average, the iterative method is a bit slower on short words but faster for long words (from a total length of 8). 
Since longer words are more numerous, using the iterative method will be advantageous.    

					\subsection{Non-commutative product}
The non-commutative and non-associative product $\vartriangleright$ is implemented as $\odot$ in \texttt{Mathematica}, an operation with no built-in meaning. We only set a few basic properties for this product:
\begin{Verbatim}[commandchars=\\£~,codes={\catcode`$=3}]
Unprotect[CircleDot]; (* $\odot$ = esc c . esc *) 
	1$\odot$x_ := x;
	x_$\odot$1 := 0;
	0$\odot$x_ := 0;
	x_$\odot$0 := 0;
Protect[CircleDot];
\end{Verbatim}
Since $\odot$ is only used as a placeholder, its actual definition in terms of shuffle product (eq. \ref{eq:productDefinition}) is coded in another function, \verb|NCP|:
\begin{Verbatim}[commandchars=\\£~,codes={\catcode`$=3}]
	NCP[0, j[b_List]] := 0;
	NCP[1, j[b_List]] := j[b];
	NCP[j[a_List], 0] := 0;
	NCP[j[a_List], 1] := j[a];
	NCP[j[a_List], j[{}]] := j[a];
	NCP[j[a_List], j[b_List]] := 
			Ap[j[a]*j[Drop[b, -1]], Last[b]] /; Length[b] > 0;
	NCP[n_*j[a_List], j[b_List]] := n*NCP[j[a], j[b]];
	NCP[j[a_List], n_*j[b_List]] := n*NCP[j[a], j[b]];
	NCP[n_*j[a_List], m_*j[b_List]] := n*m*NCP[j[a], j[b]];
	NCP[x_ + y_, z_] := NCP[x, z] + NCP[y, z];
	NCP[x_, y_ + z_] := NCP[x, y] + NCP[x, z];
	NCP[x_*(y_ + z_), t_] := NCP[x*y, t] + NCP[x*z, t];
	NCP[(y_ + z_)*x_, t_] := NCP[x*y, t] + NCP[x*z, t];
\end{Verbatim}

\subsection{Picard iteration}
The usual Picard iteration is implemented with a helper function \verb|PicardIteration| as following:
\begin{Verbatim}
PicardIteration[f_List,X_]:=
		Total[MapIndexed[(Ap[#1[X]*j[{}],First[#2]-1])&,f]]
Picard[f_List,X0_,n_Integer]:=
		Nest[(PicardIteration[f,#])&,X0,n]
\end{Verbatim}
For example, \verb|Picard[f,x0,4]| outputs the stochastic approximation of the SDE with the functions $f$ collected in a list in the first argument. This was used in \cite{Papavasiliou2010} for a system with linear drift and quadratic variance.

With the new representation in $Q$, it can be written directly as in eq. \ref{eq:picard}:
\begin{Verbatim}[commandchars=\\£~,codes={\catcode`$=3}]
PicardQ1Dim[Q_, R_, q_] :=  
		Nest[Q[0] + Sum[(#^r)$\odot$Q[r], {r, 1, q}] &, Q[0], R - 1];
\end{Verbatim}
\verb|R| is the number of iterations to be calculated and \verb|q| the maximum degree of the polynomials $f$. \verb|PicardQ1Dim[]| produces a very compact representation of the expansion, which needs to be processed further in order to give the same result as \verb|Picard[]|.

					\subsection{Distributed processing of monomials}
Going from the compact representation provided by \verb|PicardQ1Dim[]| to the linear combination of iterated integrals $j$'s is done in a few stages. Each stage modifies the representation of the stochastic expansion in such a way that a) computational requirements are minimized and b) it can be parallelized.

As described in section \ref{subsection:approach}, the workflow goes as follows:
\begin{enumerate}
\item \verb|PicardQ1Dim[]| produces a polynomial in $Q$ and $\odot$. 
\item Each monomial (in $Q$ and $\odot$) is extracted and stored in a list (actually a file).
\item For each monomial, the $Q$'s are replaced by their values in $j$ (eq.\ref{eq:QDefinition}). They are now polynomials in $j$ and $\odot$.
\item The monomials (in $j$ and $\odot$) from each polynomial are extracted and stored in a file.
\item The product $\odot$ is instanciated in terms of the shuffle product (eq. \ref{eq:productDefinition}). The result is a linear combination of iterated integrals $j$ for each monomials. The stochastic expansion is the overall sum of the all those linear combinations.
\end{enumerate}

\begin{figure}
\label{fig:workflow}
\begin{center}
\includegraphics[width=0.75\textwidth]{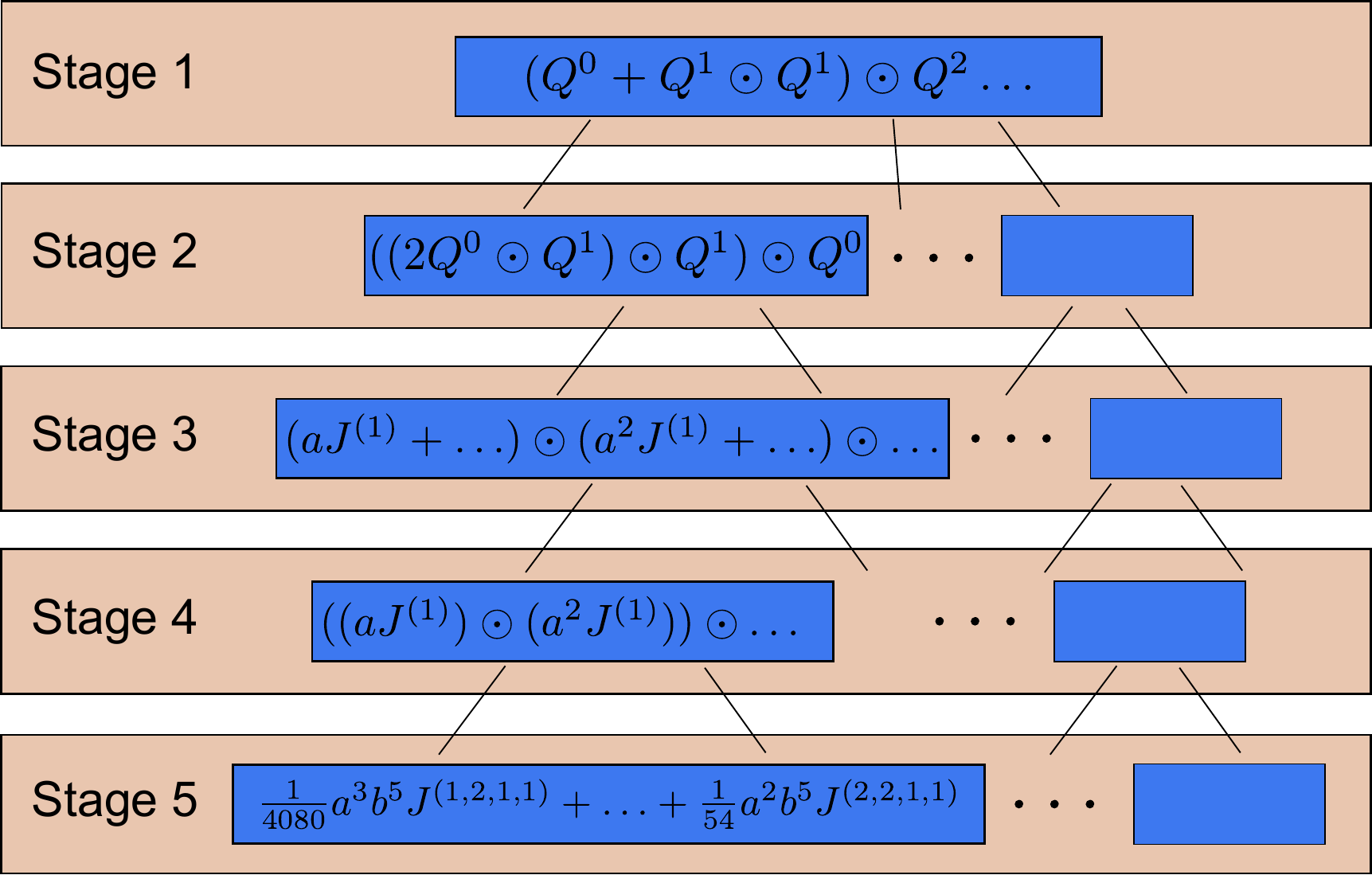}
\end{center}
\caption{Workflow. From the compact representation in $Q$ and $\odot$ to the linear combination of iterated integrals $J$'s.}
\end{figure}

As can be seen on figure \ref{fig:workflow}, the different polynomials are successively expanded in terms of monomials, which in turn are processed independently. Since each expression is effectively broken down in small parts and dumped into a file, a much larger of terms can be computed, as can be seen in the next section.

Since each processing stage follows roughly the same logic, with slightly different transformation rules, we only detail the first one: going from the polynomial in $Q$ and $\odot$ to its monomials. This is done step by step, and not by applying a global rule, since doing so would produce a very large output and run the risk of failling due to memory limits.
\begin{Verbatim}[commandchars=\\£~,codes={\catcode`|=3},numbers=left]
Print["=> Deriving the monomials in Q and |\odot|"];
changeHappened = True;
workingFile = "uexpansion";
finalFile = "uexpansion_Final";
DistributeDefinitions[workingFile, finalFile, rulePlus, rulePower,CircleDot];
(* write expansion in *)
in = OpenWrite@mainFile; Write[in, PicardQ1Dim[Q, 4, 2]];  Close@in;
SetSharedVariable[changeHappened];
ParallelEvaluate[ Put[finalFile <> ToString[$KernelID]];];
While[changeHappened,
  changeHappened = False;
  positionList = GetStreamPositions@workingFile;
  DistributeDefinitions[positionList];
  ParallelEvaluate[in = OpenRead@workingFile];
  ParallelEvaluate[Put["temporary_" <> ToString[$KernelID]]];
  WaitAll@Table[
    With[{sp = sp},
     ParallelSubmit[
      SetStreamPosition[in, sp];
      v = Read[in, Expression];
      If[FreeQ[v, Plus] && FreeQ[v, Power],
       PutAppend[v, finalFile <> ToString[$KernelID]],
       If[Head[v] === Plus,
        changeHappened = True;
        Scan[PutAppend[#, "test_temporary_" <> ToString[$KernelID]] &, v], 
        firstPlus = Min@Map[Length@# &, Position[v, Plus]];
        firstPower = Min@Map[Length@# &, Position[v, Power]];
        If[{firstPlus, firstPower} != {Infinity, Infinity},
         If[firstPlus <= firstPower,
          changeHappened = True;
          Do[v = Replace[v, rulePlus, depth], {depth, firstPlus - 2, 1, -1}];
          v = Replace[v, rulePlus];          
          Scan[PutAppend[#, "temporary_" <> ToString[$KernelID]] &, v],
          changeHappened = True;
          v = Replace[v, rulePower, {firstPower - 1}];
          PutAppend[v, "temporary_" <> ToString[$KernelID]]], 
         PutAppend[v, "temporary_" <> ToString[$KernelID]]]]]]], {sp, 
     positionList}];
  ParallelEvaluate[Close@in];
  ConcatenateFiles[tempFileNames, workingFile]];
ConcatenateFiles[Map[finalFile <> ToString@# &, kernelIDList], finalFile];
\end{Verbatim}

The data to process is in \verb|workingFile|. At each step, an entry is read from this file and processed by one of the kernels (lines 20-37). If the result is a monomial, it is added to the final file (line 22). Otherwise, it still requires further processing and is added to a temporary file (lines 25,33,36,37). When all entries in \verb|workingFile| have been processed, \verb|workingFile| is replaced by the temporary file (line 40). The algorithm loops back until all monomials have been extracted to \verb|finalFile|. Note that each kernel has its own temporary and final files, as writing to the same file in parallel typically results in corrupted files and missing entries. Processing an entry consists in finding the highest 'Plus' in the expression and distributing the monomials around it. In this manner, the polynomial is first decomposed into the largest polynomials, which are processed further separately.

\section{Expectation of an iterated integral}
It is often of interest to calculate the moments of the solution of the SDE and this can be approximated by computing the moments of the stochastic expansion. Since the expansion is a weighted sum of iterated integrals $j$, its expectation is simply the weighted sum of the integrals' expectations. 
If the drivers consist of time and Brownian motions, the expectation of an iterated integral has a simple analytic form that can be arrived at recursively (\cite{Tocino2009}). 

Here we present a more direct way of calculating this quantity. Given a word $\alpha$ and assuming Wiener processes (\cite{Ladroue2010}):
$$EJ^\alpha(t)=\left\{
\begin{array}{l}
0\textnormal{ if }\alpha \textnormal{ is not a sequence of } 0 \textnormal{ and pairs } mm\\
p_\alpha\frac{t^{q_\alpha}}{q_\alpha!}\textnormal{ otherwise}
\end{array}\right.$$

where $p_\alpha=\frac{1}{2}^\frac{\#\{\alpha_i\neq0\}}{2}$ and $q_\alpha=\frac{1}{2}(\#\{\alpha_i\neq0\})+(\#\{\alpha_i=0\})$.
\\Thus, for example:
\\$\begin{array}{lclcl}
E J^{(0,1,1,0,0)}&=&1/2^{2/2} t^{(3+2/2)}/(3+2/2)!&=&\frac{t^4}{48}\\
E J^{(0,1,1,0,0,1)}&=&0&&\\
E J^{(2,2,1,1,3,3)}&=& 1/2^{6/2} t^{(0+6/2)}/(0+6/2)!&=&\frac{t^3}{48}\\
E J^{(2,2,0,1,1,3,3,0,0,0)}&=& 1/2^{6/2} t^{(4+6/2)}/(4+6/2)!&=&\frac{t^7}{8.7!}\\
\end{array}
$

This result is implemented in \texttt{Mathematica}. The expectation for a word $\alpha$ is then calculated in at most $|\alpha|$ steps:
\begin{Verbatim}[numbers=left]
ExpSBM[t_, j[a_List]] := Module[{i, c},
   i = Length@a;
   c = {0, 0};
   Catch[
    While[i > 0,
     If[a[[i]] == 0,
      c += {0, 1}; i--,
      If[(i > 1) && (a[[i]] == a[[i - 1]]),
       c += {1, 1}; i -= 2,
       c = {Infinity, 0}; Throw@0
       ]]]];
   (1/2)^First@c t^Last@c/(Last@c)!];
\end{Verbatim}

					\section{Example}
Consider the following SDE: $dY_t=a(1-Y_t)dX_1+bY_t^2dX_2$ and $Y_0=0$. In this case, $m=1$, $n=2$, $q=2$. The two functions $f$ are $f_{1,1}(x)=a(1-x)$ and $f_{1,2}(x)=bx^2$. 

Only two things are required from the user: the definition of the objects $Q$ in a transformation rule, easily obtained by derivation (eq.\ref{eq:QDefinition}), and the number of Picard iterations to be computed.
In this case, the three $Q$'s are:
\begin{eqnarray*}
Q^0&=&\frac{1}{0!}(a(1-0)X^{(1}+b0^2X^{(2)})\\
&=& aX^{(1)}\\
Q^1&=&\frac{1}{1!}(-aX^{(1)}+2b0X^{(2)})\\
&=& -aX^{(1)}\\
Q^2&=&\frac{1}{2!}(0X^{(1}+2bX^{(2)})\\
&=& bX^{(2)}\\
\end{eqnarray*}
Therefore, the transformation rule corresponding to this system is:
\begin{verbatim}
ruleModel={Q[0]->aj[{1}],Q[1]->-aj[{1}],Q[2]->bj[{2}]}};
\end{verbatim}
If the first driver is time, we can follow the convention that drivers are numbered from 0:
\begin{verbatim}
ruleModel={Q[0]->aj[{0}],Q[1]->-aj[{0}],Q[2]->bj[{1}]}};
\end{verbatim}
A small number (\emph{e.g.} 4) of Picard iterations is sufficient for a good approximation (see \cite{Papavasiliou2010}). An optional parameter is the maximum word length of the iterated integrals.

Running the algorithm with these parameters, we obtain the expansion, which is a weighted sum of 676 iterated integrals and already has a \verb|ByteCount|  (size) of 394'288. Its expectation is a much smaller expression:
\begin{eqnarray*}
a T-\frac{a^2}{2}T^2+\frac{a^3}{6} T^3+(\frac{1}{4} a^3 b^2 -\frac{a^4}{24}) T^4-\frac{7}{20} a^4 b^2 T^5+\frac{61}{360} a^5 b^2 T^6
+(\frac{17}{140} a^5 b^4-\frac{1}{24} a^6 b^2 ) T^7&&\\
+(\frac{1}{192} a^7 b^2 -\frac{21}{160} a^6 b^4 ) T^8+\frac{157 a^7 b^4 T^9}{3024}+(\frac{43 a^7 b^6}{1800}-\frac{17 a^8 b^4}{2800})T^{10}-\frac{1}{100} a^8 b^6 T^{11}&&\\
\end{eqnarray*}

This is confirmed by the previous implementation of Picard iteration with the simple code:
\begin{Verbatim}[commandchars=\\£~,codes={\catcode`$=3}]
P0[x_] := a (1 - x);
P1[x_] := b x^2;
expansion = Picard[{P0, P1}, 0, 4];
ExpSBM[t, expansion]
\end{Verbatim}

However, if we now set the initial value to an arbitrary $y_0$ instead of 0, the stochastic expansion contains a much larger number of terms. It is calculated in the same manner using $Q$ as:
\begin{verbatim}
ruleModel={Q[0]->a(1-y0)j[{0}]+by0^2j[{1}], Q[1]->-aj[{0}]+2by0j[{1}],Q[2]->bj[{1}] };
\end{verbatim}
This results in an expansion with 10'710 unique iterated integrals and a \verb|ByteCount| of 15'841'624. This cannot be confirmed by the simple code as it runs out of memory and crashes before completing. Note that the expansion is usually left in a file whose entries are parts of the linear combination. Thus, it is possible to, for example, compute the expectation of the expansion without having to store it in memory at any point. Moreover, while computationaly expensive, these expansions can be calculated once in a general case and saved in a file for future application without having to recompute them {\it de novo}.


\section{Conclusion}
Stochastic expansions provide a local approximation of the solution of a stochastic (or deterministic) differential equation. They can be used for a variety of applications, from simulation to parameter estimation. However, as the number of terms grows exponentially with the desired precision, they can rapidly become unwiedly to manipulate. 

We presented a new way of calculating these expansions that bypasses the limitation of the usual approach, \emph{via} a re-parametrisation of the problem and the parallelization of the computation. We have shown that in a simple example our method was able to compute the expansion when a direct approach failed. We also presented two new approaches for efficiently deriving the shuffle product of two iterated integrals and the expectation of an iterated integral, when the drivers are time and Brownian motion.

So far, our approach has been implemented for one-dimensional differential equation\footnote{\texttt{Mathematica} package available on the author's webpage: \url{http://go.warwick.ac.uk/roughpaths}}. However, the theoretical foundation for the multi-dimensional case is available, as presented in section \ref{section:generalization}. Now that the computing requirements have been alleviated, an implementation for the general case is possible. Stochastic expansions will then be available for more complex systems. 

\section*{Acknowledgements}
This work was funded by the EPSRC (EP/H019588/1, 'Parameter Estimation for Rough Differential Equations with Applications to Multiscale Modelling').
 We would like to thank the \texttt{Mathematica} newsgroup (\verb|comp.soft-sys.math.mathematica group|) for their help and advice on file parallelization. The code was tested on Buster, the cluster of the department of Statistics of the university of Warwick.

\bibliographystyle{plain}
\bibliography{distributedprocedure}

\begin{thebibliography}{1}

\bibitem{Kloeden91}
P.~E. Kloeden and E.~Platen.
\newblock Relations between multiple ito and stratonovich integrals.
\newblock {\em Stochastic Analysis and Applications}, 9(3):311--321, 1991.

\bibitem{Ladroue2010}
Christophe Ladroue.
\newblock Expectation of stratonovich iterated integrals of wiener processes.
\newblock Aug 2010.

\bibitem{Lyons2003}
Terry Lyons and Zhongmin Qian.
\newblock {\em System Control and Rough Paths (Oxford Mathematical
  Monographs)}.
\newblock {Oxford University Press, USA}, February 2003.

\bibitem{Papavasiliou2010}
Anastasia Papavasiliou and Christophe Ladroue.
\newblock Parameter estimation for rough differential equations.
\newblock May 2010.

\bibitem{Tocino2009}
A.~Tocino.
\newblock Multiple stochastic integrals with mathematica.
\newblock {\em Mathematics and Computers in Simulation}, 79(5):1658--1667,
  January 2009.

\bibitem{Wolfram2003}
Stephen Wolfram.
\newblock {\em The Mathematica Book, Fifth Edition}.
\newblock Wolfram Media, 5th edition, August 2003.

\end{thebibliography}
\end{document}